\documentclass{amsart}

\newtheorem{theorem}{Theorem}[section]

\theoremstyle{definition}

\theoremstyle{remark}

\numberwithin{equation}{section}

\begin{document}

\title[On the Gaps between Zeros of the Riemann zeta-function]{A Note on the Gaps between Consecutive Zeros of the Riemann zeta-function}

\author{H. M. Bui}
\address{Mathematical Institute, University of Oxford, Oxford, OX1 3LB UK}
\email{hung.bui@maths.ox.ac.uk}
\thanks{HMB is supported by an EPSRC Postdoctoral Fellowship. MBM is supported in part by a University of Mississippi College of Liberal Arts summer research grant.  NN is supported in part by a an NSERC Discovery grant.}

\author{M. B. Milinovich}
\address{Department of Mathematics, The University of Mississippi, University, MS 38677 USA}
\email{mbmilino@olemiss.edu}

\author{N. Ng}
\address{Department of Mathematics and
Computer Science, University of Lethbridge, Lethbridge, AB Canada T1K 3M4}
\email{ng@cs.uleth.ca}





\begin{abstract}
Assuming the Riemann Hypothesis, we show that infinitely often consecutive non-trivial zeros of the Riemann zeta-function differ by at most 0.5155 times the average spacing and infinitely often they differ by at least 2.69 times the average spacing.  
\end{abstract}

\maketitle

\section{Introduction}
Let $\zeta(s)$ be the Riemann zeta-function.  Assuming the Riemann Hypothesis,  we denote the non-trivial zeros of $\zeta(s)$ as $\rho=\tfrac{1}{2}\pm i \gamma$ where $\gamma\in\mathbb{R}$.  It is well known that, for $T\geq 10$, 
\begin{equation*}\label{NT}
N(T): = \sum_{0<\gamma\leq T} 1= \frac{T}{2\pi} \log \frac{T}{2\pi} - \frac{T}{2\pi} + O\big(\log T\big). 
 \end{equation*}
Hence, if we let $\gamma \leq \gamma'$ denote consecutive ordinates of the zeros of $\zeta(s)$, we see that the average size of $\gamma'-\gamma$ is $2\pi/\log \gamma$. Normalizing, we let
$$ \lambda := \limsup_{\gamma \geq 0} \frac{(\gamma'-\gamma)\log \gamma}{2 \pi}$$
and
$$  \mu := \liminf_{\gamma \geq 0} \frac{(\gamma'-\gamma) \log \gamma}{2 \pi}$$
and we observe that  $\mu\leq 1 \leq \lambda.$  It is expected that there are arbitrarily large and arbitrarily small (normalized) gaps between consecutive zeros of the Riemann zeta-function on the critical line; in other words,  that $\mu=0$ and $\lambda=+\infty$.  
In this note, we prove the following theorem.

\begin{theorem} \label{th}
Assume the Riemann Hypothesis.  Then $\lambda>2.69$ and $\mu<0.5155.$
\end{theorem}

We briefly describe the history of the problem.  Very little is known unconditionally;  however, Selberg (unpublished, but announced in \cite{S}) has shown that $\mu<1<\lambda$.  Assuming the Riemann Hypothesis, numerous authors \cite{CGG1, H, Mo, MO, Mu} have obtained explicit bounds for $\mu$ and $\lambda$.  Theorem \ref{th} improves the previously best known results under this assumption which were $\mu<0.5172$  \ due to Conrey, Ghosh \& Gonek \cite{CGG1} and $\lambda>2.63$  \ due to R. R. Hall \cite{H}.  Assuming the generalized Riemann Hypothesis for the zeros of Dirichlet $L$-functions, Conrey, Ghosh \& Gonek \cite{CGG2} have shown that $\lambda>2.68$.  Their method can be modified (see \cite{N} and \cite{B}) to show that $\lambda>3$.

Understanding the distribution of the zeros of the zeta-function is important for a number of reasons.  One reason, in particular, is the connection between the spacing of the zeros of $\zeta(s)$ and the class number problem for imaginary quadratic fields.  This is described by Conrey \& Iwaniec in \cite{CI}; see also Montgomery \& Weinberger \cite{MW}. Studying this connection led Montgomery \cite{Mo} to investigate the pair correlation of the ordinates of the zeros of the zeta-function.  He conjectured that, for any fixed $0<\alpha<\beta$, 
$$  \sum_{\substack{0<\gamma, \gamma'\leq T \\ \frac{2\pi \alpha}{\log T}\leq\gamma'-\gamma\leq \frac{2\pi \beta}{\log T}}} \!\!\!\!\!\!\!\! 1 \
 \sim  \ N(T) \ \int_\alpha^\beta \Big( 1-\Big( \frac{\sin \pi u}{\pi u} \Big)^2 \Big)  \ du .$$
Here $\gamma$ and $\gamma'$ run over two distinct sets of ordinates of the non-trivial zeros of $\zeta(s)$.  Clearly, Montgomery's conjecture implies that $\mu=0$.  Moreover, F. J. Dyson observed that the eigenvalues of large, random  complex Hermitian or unitary matrices have the same pair correlation function.  This observation (among other things) has led to a stronger conjecture that the zeros of the zeta-function should behave, asymptotically, like the eigenvalues of large random matrices from the Gaussian Unitary Ensemble.  These ideas lead to the conjecture that $\lambda=+\infty$.

\section{Montgomery \& Odlyzko's method for exhibiting irregularity in the gaps between consecutive zeros of $\zeta(s)$}

Throughout the remainder of this note, we assume the truth of the Riemann Hypothesis. \\

 Let $T$ be large and put $K=T(\log T)^{-2}$.  Further, let
$$ h(c) := c - \frac{\text{Re } \big( \sum_{nk\leq K} a_k \overline{a_{nk}} g_c(n) \Lambda(n) n^{-1/2} \big)}{\sum_{k\leq K} |a_k|^2} $$
where 
$$ g_c(n) = \frac{2 \sin\big( \pi c \frac{\log n}{\log T}\big)}{\pi \log n}$$
and $\Lambda(\cdot)$ is von Mangoldt's function defined by $\Lambda(n)=\log p$ if $n=p^k$ for a prime $p$ and by $\Lambda(n)=0$, otherwise.  In \cite{MO}, by an argument using the Guinand-Weil explicit formula for the zeros of $\zeta(s)$, 
Montgomery \& Odlyzko show that if $h(c)<1$ for some choice of $c>0$ and a sequence $\{a_n\}$ then $\lambda \geq c$ and that if  $h(c)>1$ for a choice of $c>0$ and a sequence $\{a_n\}$ then $\mu \leq c$.  In particular, for any such choices of $c$ and $\{a_n\}$, their method proves the existence of a pair of consecutive zeros of $\zeta(s)$ with ordinates $\gamma\leq \gamma'$ in the interval $[T/2, 2T]$ which satisfy $\gamma'-\gamma \geq \frac{2\pi c}{\log T}$ and $\gamma'-\gamma \leq \frac{2\pi c}{\log T}$, respectively.

Conrey, Ghosh \& Gonek \cite{CGG1} have given an alternative, and much simpler, way of viewing this problem.\footnote{Mueller \cite{Mu} was the first to observe that this method gave a lower bound for $\lambda$. In \cite{CGG1}, it was noticed that the method also gave an upper bound for $\mu$ and was equivalent to the method in \cite{MO}.}  Let
$$ A(t) = \sum_{k\leq K} a_k k^{-it}$$
be a Dirichlet polynomial and set 
$$ M_1= \int_T^{2T} \big|A(t)\big|^2 \ \! dt  \quad \text{ and } \quad M_2(c) = \int_{-\pi c/\log T}^{\pi c/\log T} \sum_{T\leq \gamma \leq 2T} \big|A(\gamma+\alpha)\big|^2 \ \! d\alpha.$$
Then, clearly, $M_2(c)$ is monotonically increasing and $M_2(\mu)\leq M_1 \leq M_2(\lambda)$.  Therefore, if it can be shown that $M_2(c)<M_1$ for some choice of $A(t)$ and $c$, then $\lambda>c$.  Similarly, if $M_2(c)>M_1$ for some choice of $A(t)$ and $c$, then $\mu<c$.  Using standard techniques to estimate $M_1$ and $M_2(c)$, it can be shown that $$M_2(c)/M_1 = h(c) +o(1).$$
Hence, this argument is seen to be equivalent to Montgomery \& Odlyzko's method, described above.  Moreover, we note that this formulation of the method suggests that we should choose a test function $A(t)$ which is small near the zeros of $\zeta(s)$ to exhibit large gaps between the zeros of the zeta-function and a test function $A(t)$ which is large near the zeros of $\zeta(s)$ to exhibit small gaps.


In \cite{MO}, Montgomery \& Odlyzko make the choices of 
$$ a_k=\frac{1}{\sqrt{k}} f\big(\tfrac{\log k}{\log K}\big) \quad  \text{ and } \quad a_k=\frac{\lambda(k)}{\sqrt{k}} f\big(\tfrac{\log k}{\log K}\big) $$
where $f$ is a continuous function of bounded variation on $[0,1]$ normalized so that $\int_0^1 |f|^2 =1$ and $\lambda(k)$, the Liouville function, equals $(-1)^{\Omega(k)}$; here, $\Omega(k)$ denotes the total number of primes dividing $k$.  By optimizing over such functions $f$, the values 
 $\mu<0.5179$ and $\lambda>1.97$ are obtained.  The authors of \cite{MO} choose $f$ to be a certain modified Bessel function and they mention this is close
 to an optimal choice. 

In \cite{CGG1}, Conrey, Ghosh \& Gonek choose the coefficients $$ a_k=\frac{d_r(k)}{\sqrt{k}}  \quad  \text{ and } \quad a_k=\frac{\lambda(k) d_r(k)}{\sqrt{k}}$$
where $d_r(k)$ is a multiplicative function defined on integral powers of a  prime $p$ by  
$$ d_r(p^k) = \frac{\Gamma(k+r)}{\Gamma(r) k!}.$$  In this context, this becomes an optimization in the variable $r$. 
The choice $r=1.1$ yields $\mu<0.5172$ and the choice $r=2.2$ yields $\lambda>2.337$.

In order to prove Theorem \ref{th}, we combine the approaches of \cite{MO} and \cite{CGG1}. We 
choose the coefficients 
$$ a_k=\frac{d_r(k)}{\sqrt{k}} f\big(\tfrac{\log K/k}{\log K}\big)  \quad  \text{ and } \quad a_k=\frac{\lambda(k) d_r(k)}{\sqrt{k}} f\big(\tfrac{\log K/k}{\log K}\big)$$
for sufficiently smooth functions $f$.   This variant allows us to optimize over both $r$ and $f$, rather than over just $r$ or $f$. 

We now provide further insight into the choice of these coefficients.  For simplicity, suppose $f$ is a polynomial. 
 Since, for Re $s>1$, 
$$ \sum_{k\geq 1} \frac{d_r(k)}{k^s} = \zeta(s)^r \quad \text{ and } \quad \sum_{k\geq 1} \frac{\lambda(k) d_r(k)}{k^s} = \Big(\frac{\zeta(2s)}{\zeta(s)}\Big)^r,$$
with our choice of coefficients we see that the test function $A(t)$ approximates $\zeta(\tfrac{1}{2}+it)^r$ and $\zeta(1+2it)^r/\zeta(\tfrac{1}{2}+it)^r$, respectively, and should have the desired effect of making $A(t)$ small (respectively large) near the zeros of $\zeta(s)$.  
Moreover, when we multiply $d_r(k)$ by $f\big(\tfrac{\log K/k}{\log K}\big) $ then $A(t)$ behaves like a linear combination of $\zeta(\tfrac{1}{2}+it)^r$ and its derivatives
and an analogous comment applies to the other case. 
The presence of the function $f$ leads to improved numerical results for bounds for $\mu$ and $\lambda$.  

\section{Large Gaps: a lower bound for $\lambda$}

In this section, we establish a lower bound for $\lambda$ by evaluating $h(c)$ with the coefficients
 $$ a_k=\frac{d_r(k)}{\sqrt{k}} f\big(\tfrac{\log K/k}{\log K}\big) $$ 
where $f$ is a continuous, real-valued function of bounded variation on $L^2[0,1]$. In what follows, we assume that $r\geq 1$ so that $d_r(mn) \leq d_r(m) d_r(n)$ for $m,n\in \mathbb{N}$.  It is well known that, for fixed $r\geq 1$,
\begin{equation}
 \label{eq:divsum}
  \sum_{k\leq x} \frac{d_r(k)^2}{k} = A_r (\log x)^{r^2} + O\big( (\log T)^{r^2-1}\big)
\end{equation}
uniformly for $x\leq T$; here $A_r$ is a certain arithmetical constant (the exact value is not important in our argument).  By partial summation, recalling that $K=T(\log T)^{-2}$,  we find that the denominator in the ratio of sums in the definition of $h(c)$ is
\begin{equation*}
\begin{split}
\sum_{k\leq K} |a_k|^2 &= \int_{1^-}^K  f\big(\tfrac{\log K/x}{\log K}\big)^2 d\Big(  \sum_{k\leq x} \frac{d_r(k)^2}{k} \Big)
\\
&= A_r r^2 \int_1^K f\big(\tfrac{\log K/x}{\log K}\big)^2 (\log x)^{r^2-1} \frac{ \ dx}{x} + O_{f,r} \big( (\log T)^{r^2-1}\big)
\end{split}
\end{equation*}
by \eqref{eq:divsum}.  By the variable change $u = 1- \frac{\log x}{\log K}$ 
\begin{equation*}
\begin{split}
\sum_{k\leq K} |a_k|^2 &=  A_r r^2 (\log K)^{r^2} \int_0^1(1\!-\!u)^{r^2-1} f(u)^2  \ \!du + O_{f,r} \big( (\log T)^{r^2-1}\big)
\\
&= A_r r^2 (\log T)^{r^2} \int_0^1 (1\!-\!u)^{r^2-1} f(u)^2  \ \! du + O_{f,r,\varepsilon} \big( (\log T)^{r^2-1+\varepsilon}\big)
\end{split}
\end{equation*}
where $\varepsilon>0$ is arbitrary.  

We now evaluate the numerator in the ratio of sums in the definition of $h(c)$.  If we let
$$ N(c):= \sum_{nk\leq K} a_k \overline{a_{nk}} g_c(n) \Lambda(n) n^{-1/2}$$
then a straightforward argument shows that
\begin{equation*}
\begin{split}
N(c) &=\frac{2}{\pi}\sum_{nk\leq K} \frac{d_r(k)d_r(kn)\Lambda(n)}{kn \log n} f\big(\tfrac{\log K/k}{\log K}\big) f(\tfrac{\log K/nk}{\log K}\big)\sin\big(\pi c \tfrac{\log n}{\log T}\big)
\\
&= \frac{2}{\pi}\sum_{pk\leq K} \frac{d_r(k)d_r(kp)}{kp} f\big(\tfrac{\log K/k}{\log K}\big) f(\tfrac{\log K/pk}{\log K}\big)\sin\big(\pi c \tfrac{\log p}{\log T}\big) + O_{f,r} \big( (\log T)^{r^2-1}\big)
\\
&= \frac{2 r}{\pi} \sum_{p\leq K} \frac{\sin\big(\pi c \tfrac{\log p}{\log T}\big)}{p} \sum_{k\leq K/p} \frac{d_r(k)^2}{k} f\big(\tfrac{\log K/k}{\log K}\big) f(\tfrac{\log K/pk}{\log K}\big) + O_{f,r} \big( (\log T)^{r^2-1}\big).
\end{split}
\end{equation*}
By Stieltjes integration and a variable change, 
the inner sum in the main term of the last expression for $N(c)$ is 
\begin{equation*}
\begin{split}
& \int_{1^-}^{K/p}  f\big(\tfrac{\log K/x}{\log K}\big) f(\tfrac{\log K/px}{\log K}\big) d\Big(  \sum_{k\leq x} \frac{d_r(k)^2}{k} \Big)
\\
&\quad \quad=  A_r r^2  \int_1^{K/p}  f\big(\tfrac{\log K/x}{\log K}\big) f(\tfrac{\log K/px}{\log K}\big) (\log x)^{r^2-1} \frac{ \ dx}{x}
 \ \!+ \ \!O_{f,r}\big( (\log T)^{r^2-1}\big)
\\
&\quad \quad =  A_r r^2  (\log K)^{r^2}\!\! \int_{\frac{\log p}{\log K}}^{1} (1\!-\! u)^{r^2-1} f(u) f(u\!-\! \tfrac{\log p}{\log K}\big)  \ \!du \ \!
 + \ \! O_{f,r}\big( (\log T)^{r^2-1}\big).
\end{split}
\end{equation*}
By combining the above estimates and interchanging the order of summation and integration, we conclude that $N(c)=M(c)+O_{f,r}\big( (\log T)^{r^2-1}\big)$ where
\begin{equation*}
\begin{split}
M(c)&= \frac{2 A_r r^3}{\pi} (\log K)^{r^2}\int_{\frac{\log 2}{\log K}}^1 (1\!-\!u)^{r^2-1} f(u) \sum_{2 \leq p \leq K^u} \frac{\sin\big(\pi c \tfrac{\log p}{\log T}\big)}{p} f\big(u\!-\! \tfrac{\log p}{\log K}\big)  \ \! du
\\
&= \frac{2 A_r r^3}{\pi} (\log K)^{r^2}\int_{0}^1 (1\!-\!u)^{r^2-1} f(u) \sum_{2\leq p \leq K^u} \frac{\sin\big(\pi c \tfrac{\log p}{\log T}\big)}{p} f\big(u\!-\! \tfrac{\log p}{\log K}\big)  \ \! du 
\\
&\quad \quad  + O_{f,r,\varepsilon}\big((\log T)^{r^2-1+\varepsilon}\big).
\end{split}
\end{equation*}
By the prime number theorem with remainder term, it follows that
$$
\sum_{2\leq p \leq K^u} \!\!\!\frac{\sin\big(\pi c \tfrac{\log p}{\log T}\big)}{p} f\big(u\!-\!\tfrac{\log p}{\log K}\big) 
= \int_{2}^{K^u}  \frac{\sin\big(\pi c \tfrac{\log x}{\log T}\big)}{x \log x} f\big(u\!-\! \tfrac{\log x}{\log K}\big) \ \! dx+ O_{f,r}\Big(\frac{1}{\log T}\Big).
$$
By the variable change $v= \frac{\log x}{\log K}$, the integral is
$$
\int_{\frac{\log 2}{\log K}}^u \frac{\sin\big(\pi c v \tfrac{\log K}{\log T}\big)}{v} f(u\!-\!v) \ \! dv = \int_0^u \frac{\sin(\pi c v )}{v} f(u\!-\!v) dv + O_{f,r,\varepsilon}\big((\log T)^{-1+\varepsilon}\big).
$$
Hence,
\begin{equation*}
\begin{split}
M(c) &= \frac{2 A_r r^3}{\pi} (\log K)^{r^2}\int_{0}^1 (1\!-\!u)^{r^2-1} f(u)  \int_0^u \frac{\sin(\pi c v )}{v} f(u\!-\!v) \ \! dv \ \!  du
\\
& \quad \quad  + O_{f,r,\varepsilon}\big((\log T)^{r^2-1+\varepsilon}\big).
\\
\end{split}
\end{equation*}
Consequently, we find that
\begin{equation} \label{eq:opt1}
\begin{split}
 h(c) &= c - \frac{2r}{\pi} \frac{\int_{0}^1 (1\!-\!u)^{r^2-1} f(u)  \int_0^u \frac{\sin(\pi c v )}{v} f(u\!-\! v)  \ \! dv \ \! du}{ \int_0^1 (1\!-\!u)^{r^2-1} f(u)^2  \ \!du} 
 \\
 &\quad \quad + O_{f,r,\varepsilon}\big((\log T)^{-1+\varepsilon}\big).
 \end{split}
 \end{equation}
Choosing $r=3.1$ and $f(x)=1+10x+39x^2$, we obtain (by a numerical calculation) that $ h(2.69)<1$ when $T$ is sufficiently large.
This provides the lower bound for $\lambda$ in Theorem \ref{th}.

\section{Small Gaps: an upper bound for $\mu$}

Since $\lambda(n)^2=1$ and $\lambda(pn)=-\lambda(n)$ for every $n \in \mathbb{N}$ and every prime $p$, we can evaluate $h(c)$ using the coefficients
 $$ a_k=\frac{\lambda(k) d_r(k)}{\sqrt{k}} f\big(\tfrac{\log K/k}{\log K}\big) $$ 
  in a similar manner to the calculations in the previous section.   Here, as above, $f\in L^2[0,1]$ is a real-valued, continuous function of bounded variation. With this choice of coefficients, we obtain that 
\begin{equation} \label{eq:opt2}
\begin{split}
h(c) &= c + \frac{2r}{\pi} \frac{\int_{0}^1 (1\!-\!u)^{r^2-1} f(u)  \int_0^u \frac{\sin(\pi c v )}{v} f(u\!-\! v)  \ \! dv \ \! du}{ \int_0^1 (1\!-\!u)^{r^2-1} f(u)^2  \ \!du} 
\\
& \quad \quad + O_{f,r,\varepsilon}\big((\log T)^{-1+\varepsilon}\big).
\end{split}
\end{equation}
Choosing $r=1.23$ and $f(x)=1+0.99x-0.42x^2$, a numerical calculation implies that $ h(0.5155)>1.$ This provides the upper bound for $\mu$ in Theorem \ref{th}.

\section{Some Concluding Remarks}

Theorem \ref{th} offers the best known bounds for $\lambda$ and $\mu$ assuming the Riemann Hypothesis; however, we are still far from proving the conjectured values of $\mu=0$ and $\lambda=\infty$.  In fact, it known that this is not attainable using Montgomery and Odlyzko's method.  Specifically, in \cite{CGG1}, it is shown that $h(c)<1$ if $c<\tfrac{1}{2}$ and $h(c)>1$ if $c\geq 6.2$.  It would be interesting to better understand the limitations of this method and, in particular, if it can be used to show that $\mu\leq \tfrac{1}{2}$.

We have not been able to prove that our bounds for $\lambda$ and $\mu$ in Theorem \ref{th} are the optimal bounds for our choice of coefficients  
$$ a_k=\frac{d_r(k)}{\sqrt{k}} f\big(\tfrac{\log K/k}{\log K}\big)  \quad  \text{ and } \quad a_k=\frac{\lambda(k) d_r(k)}{\sqrt{k}} f\big(\tfrac{\log K/k}{\log K}\big).$$
In the special case of $r=1$, 
this optimization problem has been solved (in terms of prolate spheroidal wave functions).  When $r\neq1$, the analogous optimization problem seems considerably more difficult.  Instead of trying to solve it explicitly, we have instead chosen $f$ to be a polynomial of low degree ($\leq 6$). Using Mathematica, we numerically optimize \eqref{eq:opt1} and \eqref{eq:opt2} over
the coefficients for each choice of $c$ and $r$. 


 Our numerical calculations seem to suggest that degree 2 polynomials work well; it does not seem like there is much to gain by increasing the degree of $f$. To demonstrate this phenomenon, we observe that one can recover the bounds for $\lambda$ and $\mu$, in the case of $r=1$, derived in \cite{MO} using degree 2 polynomials in place of prolate spheroidal wave functions (which are much more difficult to compute numerically).  In this case, after some rearranging and a change of variables, our above calculations imply that 
$$  h(c) = c \pm  \frac{\int_{0}^1 f(u) \int_0^1 f(v) \frac{\sin(\pi c (u - v ))}{\pi (u - v)}  \ \! dv \ \!  du}{ \int_0^1 f(u)^2  \ \!du} + o(1).$$ 
Letting $f(x)=1+0.46526x-0.46526x^2$, a numerical calculation shows that $h(0.5179)>1$ and if we let $f(x)=1+17.9426x-17.9426x^2$, then it can be shown that $h(1.97)<1$.  These are essentially the (optimal) values obtained by Mongomery and Odlyzko in \cite{MO} when $r=1$.

\bibliographystyle{amsplain}

\begin{thebibliography}{10}

\bibitem{B} H. M. Bui, \textit{Large gaps between consecutive zeros of the Riemann zeta-function}, submitted.  Available on the \texttt{arXiv} at \texttt{http://arxiv.org/abs/0903.4007}.

\bibitem{CGG1} J. B. Conrey, A. Ghosh, \and S. M. Gonek, \textit{A note on gaps between zeros of the zeta function}, Bull. London Math. Soc. \textbf{16} (1984), 421--424.

\bibitem{CGG2} J. B. Conrey, A. Ghosh, \and S. M. Gonek,  \textit{Large gaps between zeros of the zeta-function}, Mathematika \textbf{33} (1986), 212--238.

\bibitem{CI} J. B. Conrey \and H. Iwaniec, \textit{Spacing of zeros of Hecke $L$-functions and the class number problem},  Acta Arith. \textbf{103}  (2002),  no. 3, 259--312.

\bibitem{H} R. R. Hall, \textit{A new unconditional result about large spaces between between zeta zeros}, Mathematika \textbf{52} (2005), 101--113.

\bibitem{Mo} H. L. Montgomery, \textit{The pair correlation of the zeros of the zeta function}, Proc. Symp. Pure Math. 24, A.M.S., Providence 1973, 181--193.

\bibitem{MO} H. L. Montgomery \and A. M. Odlyzko, \textit{Gaps between zeros of the zeta function}, Coll. Math. Soc. J\~{a}nos Bolyai 34.  Topics in Classical Number Theory, Budapest, 1981.

\bibitem{MW} H. L. Montgomery \and P. J. Weinberger, \textit{Notes on small class numbers}, Acta Arith. \textbf{24} (1974), 529--542.

\bibitem{Mu} J. Mueller, \textit{On the difference between consecutive zeros of the Riemann zeta-function}, J. Number Theory \textbf{14} (1983), 393--396.

\bibitem{N} N. Ng, \textit{Large gaps between the zeros of the Riemann zeta function}, J. Number Theory \textbf{128} (2008), 509--556.

\bibitem{S} A. Selberg, \textit{The zeta-function and the Riemann Hypothesis}, Skandinaviske Matematikerkongres \textbf{10} (1946), 187--200.



\end{thebibliography}

\end{document}